\documentclass{article}
\usepackage{amsmath}
\usepackage[left=3cm, right=3cm, bottom=2.8cm, top=2.8cm]{geometry}
\usepackage{hyperref}
\usepackage{cleveref}
\usepackage{amssymb}
\usepackage{mathtools}
\usepackage{tikz-cd}
\usepackage{amsthm}
\usepackage{graphicx}
\usepackage{graphics}
\usepackage{mathrsfs}
\usepackage{cleveref}
\usepackage{thmtools}
\usepackage{enumitem}
\usepackage{bbm}
\newlist{steps}{enumerate}{1}
\usepackage{mathtools}
\makeatletter
\newcommand{\xMapsto}[2][]{\ext@arrow 0599{\Mapstofill@}{#1}{#2}}
\def\Mapstofill@{\arrowfill@{\Mapstochar\Relbar}\Relbar\Rightarrow}
\makeatother
\usepackage{array}
\setlist[steps, 1]{label = Step \arabic*:}
\theoremstyle{plain}
\usepackage{footnote}
\makesavenoteenv{tabular}
\newtheorem*{theorem*}{Theorem}
\newtheorem{theorem}{Theorem}[subsection]
\newtheorem{definition}[theorem]{Definition}

\newtheorem{proposition}[theorem]{Proposition}
\newtheorem{remark}[theorem]{Remark}

\newtheorem*{corollary*}{Corollary}
\newtheorem{exmp}[theorem]{Example}

\newtheorem{lemma}[theorem]{Lemma}
\newtheorem*{lemma*}{Lemma}

\numberwithin{equation}{subsection}
\usepackage{rotating}
\usepackage{comment}
\usepackage[nottoc,notlot,notlof]{tocbibind} 

\usepackage[toc,page]{appendix}

\title{A Mathematical Definition of Path Integrals on Symplectic Manifolds} 
\author{Joshua Lackman\footnote{josh@pku.edu.cn} }
\date{}

\begin{document}

\maketitle
\begin{abstract}
\noindent We give a mathematical definition of some path integrals, emphasizing those relevant to the quantization of symplectic manifolds (and more generally, Poisson manifolds) — in particular, the coherent state path integral. We show that K\"{a}hler manifolds provide many computable examples and we emphasize those whose Bergman kernel is constant along the diagonal.
\end{abstract}
\tableofcontents
\section{Introduction}
We discuss a mathematical definition of some path integrals, particularly those relevant to the quantization of symplectic manifolds $(M,\omega),$ such as 
\begin{align}
&\label{pcon}\langle m_1|m_0\rangle=\int_{\gamma(0)=m_0}^{\gamma(1)=m_1}\mathcal{D}\gamma\,e^{\frac{i}{\hbar}\int_0^1\gamma^*\nabla}\;,
\\\label{def}&(f\star g)(m)=\int_{X(\infty)=m}\mathcal{D}X\,e^{\frac{i}{\hbar}\int_{D^2}X^*\omega}f(X(1))g(X(0))\;,
\end{align}
for $m_0,m_1,m\in M.$ In \cref{pcon}, $\nabla$ is the connection on the prequantum line bundle and the integrand should be interpreted as parallel transport over $\gamma$ from the vector space over $m_0$ to the vector space over $m_1.$\footnote{Implicit in this notation is a choice of orthonormal bases of the vector spaces over $m_0, m_1.$ We have written it this way to be consistent with the physics literature, eg. equation (2.7) of \cite{brane}. See \cref{mg} for the definition.}  Computing this path integral \textit{immediately} produces a quantization, ie. a good Hilbert space and quantization map defined on all classical observables, as explained in \cref{simple}, \cref{mg}. We will call it a propagator.
\\\\We suggest that computing \cref{pcon} is equivalent to finding a propagator satisfying \cref{prop}. In \cref{def}, the integral is over maps from the disk into $M,$ where $0,1,\infty$ are points on the boundary of the disk. It is related to \cref{pcon}, since
\begin{equation}
  e^{\frac{i}{\hbar}\int_{D^2}X^*\omega}= e^{\frac{i}{\hbar}\int_{\partial D}X\vert_{\partial D}^*\nabla}\;.
\end{equation} 
\Cref{def} naturally generalizes to Poisson manifolds. K\"{a}hler manifolds produce many computable examples, with particularly nice ones arising from those whose Bergman kernel is constant along the diagonal. As such this is closely related to Toeplitz quantization \cite{bord}, as discussed in \cite{charles}, \cite{sergio}.
\\\\\Cref{pcon} is sometimes called the coherent state path integral (or phase space path integral) and Daubechies, Klauder in \cite{klauder}, \cite{klauder2} computed it on surfaces of constant curvature using Brownian motion — see \cref{surf}, where we also discuss the issue of uniqueness.\footnote{It was more generally computed on K\"{a}hler manifolds in \cite{charles} — when we say it was computed, we mean a reasonable value was assigned to it. A general definition isn't given.} In some sense it is the most natural and mathematically well-behaved form of quantization, and as we will discuss, it is closely related to the different geometric and deformation quantization programs, \cite{woodhouse}, \cite{weinstein}, \cite{kontsevich}, \cite{fedosov}. In particular, \cref{def} is closely related to the Poisson sigma model, see \cite{catt} remark 3, and \cite{bon} equation (3.39). Also see \cite{grady}.
\\\\The idea is that we can compute \cref{pcon} by finding a function on $M\times M$ with the right derivative.  We will not prove any very general existence theorems, but we do obtain well-defined existence and uniqueness questions. This paper is largely a continuation of \cite{Lackman2}, where we formalized the procedure of putting path integrals over morphisms of Lie algebroids on a lattice (eg. path integrals over maps between manifolds, \cref{lat}). That paper describes a path integral construction analogous to the construction of the Riemann integral via Riemann sums, eg. Feynman's construction of the path integral. This paper, on the other hand, focuses on the path integral analogue of the fundamental theorem of calculus — the construction of the propagator fits into the framework of the aforementioned paper.
\\\\In particular, we have the following lemma:
\begin{lemma*}(\ref{berg})
The normalized Bergman kernel $B$ of the prequantization of a K\"{a}hler manifold integrates the prequantum connection (ie. exactly equals \cref{pcon}) if and only if $B$ is constant along the diagonal.
\end{lemma*}
Indeed, the examples computed by Daubechies, Klauder all have the property that the Bergman kernel is constant along the diagonal. The significance in geometric quantization of it being constant has been known at least since \cite{cahen} (where they also give a sufficient condition for it to be constant), and such prequantizations were studied by Donaldson in \cite{donald}, where they are called balanced. In general, the normalized Bergman kernel integrates the prequantum connection of the Bergman K\"{a}hler form, which is an $\mathcal{O}(\hbar)$-perturbation of the original K\"{a}hler form. It is a rescaling of the pullback of the Fubini-Study form by the Kodaira embedding, \cite{lu1}, \cite{ruan}.
\\\\$\textbf{Idea:}$
\\\\The definition we give is simple, and to understand it we will begin with an analogy with the definition of antiderivatives, which will lead us into \cref{condition1}, \ref{condition2}:
\\\\Let $f\,dx$ be a 1-form on $\mathbb{R}.$ An antiderivative is normally said to be a function $F$ such that $dF=f\,dx.$ However, an antiderivative can be equivalently defined as a function 
\begin{equation}
    F:\mathbb{R}\times\mathbb{R}\to\mathbb{R}
\end{equation}
which satisfies the following two conditions:
\begin{align}\label{cond}
    & F(x,z)+F(z,y)=F(x,y)\,,
    \\\label{cond2}& \partial_y F(x,y)\vert_{y=x}=f(x)\,.
\end{align}
The fundamental theorem of calculus says that $F(x,y)=\int_x^y df\,.$ 
\\\\Suppose that we want to find such an antiderivative. In order to compute it, we can first choose \textit{any} function $F$ satisfying \cref{cond2}, ignoring \cref{cond}. We can then triangulate the interval $[x,y]$ with vertices $x=x_0<\cdots<x_n=y$ and consider the sum
\begin{equation}
    \sum_{i=0}^{n-1}F(x_i,x_{i+1})\;.
\end{equation}
This gives an approximation to $\int_x^y df\,,$ whose limit over subdivisions converges to $\int_x^y df\,.$ The traditional choices for $F$ are the ones giving the left-point and right--point Riemann sums:
\begin{align}
    & F(x,y)=f(x)(y-x)\,,
    \\& F(x,y)=f(y)(y-x)\,.
\end{align}
However, any choice of smooth $F$ satisfying \cref{cond2} will do, see \cref{int}. This gives a notion of Riemann sum which behaves well with respect to pullbacks, which is needed for the path integral. We can think of a path integral as an antiderivative of a differential form in a different category — a category of complex measures. Feynman's way of computing them is the analogue of taking a limit of Riemann sums.
\\\\Note that, we can replace the additive structure used in \cref{cond} with a multiplicative structure:
\begin{align}
    &\label{condm1} F(x,z)F(z,y)=F(x,y)\,,
    \\& \partial_y\log{F(x,y)}\vert_{y=x}=f(x)\;,
    \\& F(x,x)=1\;.
\end{align}
However, this doesn't really give anything different, but we do get something different if we only ask that \cref{condm1} holds ``on average". This leads us to the path integral.
\begin{remark}
 Stated differently, an antiderivative of $f\,dx$ is a 1-cocycle $F$ on the pair groupoid 
\begin{equation}
    \textup{Pair}\,\mathbb{R}\rightrightarrows\mathbb{R}
    \end{equation}
such that $VE_0(F)=f\,dx\,,$ where $VE_0$ is the van Est map. See appendix \ref{Riemann}.  
\end{remark}
\section{The Simplest Case}\label{simple}
\subsection{Definition of the Path Integral}
Let $(M,d\mu)$ be a manifold with a measure $d\mu$ and consider a purely imaginary 1-form $i\omega.$ Let 
\begin{equation}
    \Omega:M\times M\to \mathbb{C}\,,
\end{equation} 
be a smooth function, considered as an integral kernel with respect to $d\mu.$ Suppose that $\Omega$ satisfies the following three conditions (compare with \cref{cond}, \ref{cond2}): for all $(x,y)\in M\times M\,,$
\begin{align}
   \label{condition1} &\int_M \Omega(x,z)\Omega(z,y)\,d\mu(z)=\Omega(x,y)\,,
    \\\label{condition2}& d_y\log{\Omega}(x,y)\vert_{y=x}=i\omega\vert_x\,,
    \\&\Omega(x,x)=1\;,
    \\\label{cond55}& \overline{\Omega(x,y)}=\Omega(y,x)\,.
\end{align}
The first condition says that the convolution $\Omega\ast \Omega=\Omega,$ and this is one of the defining properties of a (time-independent) propagator;\footnote{See \cref{brown} for a comment on the time-dependent case.} $d_y$ is the exterior derivative in the second factor; the additional third and fourth conditions are normalization and reversal of orientation properties.\footnote{These condition is automatically satisfied by $iF\,,$ where $F$ satisfies \cref{cond} and with the multiplicative structure replaced by the additive structure.}
\\\\
The first property can be thought of as a cocycle condition on the pair groupoid $\textup{Pair}\,M\,,$ and the second condition says that $VE_0(\log{\Omega})=i\omega\,,$ where $VE_0$ is the van Est map, see appendix \ref{Riemann} (note that, $\Omega$ is non-zero on a neighborhood of the identity bisection).
\begin{definition}\label{prop1}
If $\Omega$ satisfies the previous four conditions, then we say that $\Omega$ integrates $i\omega$ (as a propagator), or that $\Omega$ differentiates to $i\omega\,.$     
\end{definition} $\Omega(x,y)$ should be thought of as the amplitude of a particle beginning at $x$ to end at $y.$\footnote{With a zero Hamiltonian.} To understand in what sense $\Omega$ is a path integral, suppose that 
\begin{equation}
    \Omega(x,y)=e^{S(x,y)}\,.
    \end{equation}
Using \cref{condition1}, we have that
\begin{equation}
    \Omega(x,y)=\underbrace{(\Omega\ast\cdots\ast\Omega)}_{n \text{ times}}(x,y)\,,
\end{equation}
hence
\begin{equation}\label{approx}
    e^{S(x,y)}=\int_{M^{n-1}}e^{\sum_{k=0}^{n-1} S(x_k,x_{k+1})}\,d\mu(x_1)\,\ldots\,d\mu(x_{n-1})\;.
\end{equation}
Now pick a smooth path $\gamma:[0,1]\to M$ and pick a triangulation $0=t_0<\ldots<t_n=1$ of $[0,1]\,.$ Due to \cref{condition2} (see \cref{int}), we have that 
\begin{equation}
    \sum_{k=0}^{n-1} S(\gamma(t_k),\gamma(t_{k+1}))\xrightarrow[]{\Delta t_k\to 0}i\int_0^1\gamma^*\omega\;.
\end{equation}
Hence, taking $n\to\infty$ in \cref{approx}, we can formally say that
\begin{equation}\label{ome}
\Omega(x,y)=\int_{\gamma(0)=x}^{\gamma(1)=y}e^{i\int_0^1\gamma^*\omega}\,\mathcal{D}\gamma\;.
\end{equation}
Now, if we ignore condition \cref{condition1}, then the limit of the right side of \cref{approx} is essentially Feynman's method of defining \cref{ome}, which is made precise by the lattice framework developped in \cite{Lackman2}, see \cref{lat}. Formally then, the right side of \cref{approx} is a Riemann sum with the algebraic property of \cref{cond} replaced by \cref{condition1} (of course, we have to understand Riemann sums more generally to make this statement, see \cref{new}).
\\\\By the complex version of Komolgorov's extension theorem (A.2 of \cite{long}), such an $\Omega$ defines a complex measure on the space of paths $[0,1]\to M\,,$ though a priori this measure is concentrated on all paths (including discontinuous ones).\footnote{One should be able to modify the measurable space so that this measure is concentrated on a smaller set of discontinuous paths. A relevant discussion is in Berezin \cite{ber}. In the case of a non-trivial line bundle, as in the next section, the measure is operator-valued.} We can denote this measure by 
\begin{equation}
    e^{i\int_0^1\gamma^*\omega}\,\mathcal{D}\gamma\;.
\end{equation}
Most of the observables in quantum theory are given by correlation functions and are of the form
\begin{equation}
    \int_{\gamma(0)=x}^{\gamma(1)=y}f_1(\gamma(t_1))\ldots f_n(\gamma(t_n))\,e^{i\int_0^1\gamma^*\omega}\,\mathcal{D}\gamma\;,
\end{equation}
for some functions $f_1,\ldots,f_n:M\to\mathbb{C}$ and $t_1\le\cdots\le t_n\in (0,1)\,.$ This is equal to
\begin{equation}
 \int_{M^{n}}f_1(x_1)\cdots f_n(x_n)\,\Omega(x,x_{1})\cdots\Omega(x_n,y)\,d\mu(x_1)\,\ldots\,d\mu(x_{n})\;.\footnote{Presently, we are only considering time-independent propagators (or equal-time propagators) hence the values of $t_1,\ldots, t_n$ don't affect the result.}   
\end{equation}
\begin{remark}\label{brown}
In the time-dependent case, eg. for Brownian motion or a non-zero Hamiltonian, \cref{condition1} is replaced by
\begin{equation}
    \int_M \Omega_{t_0,t_1}(x,z)\Omega_{t_1,t_2}(z,y)\,\omega_z^n=\Omega_{t_0,t_2}(x,y)\;.
\end{equation}

\end{remark}
\subsection{The Hilbert space and Path Integral – Operator Correspondence}\label{corr}
This section puts the work of \cite{klauder2} into general context. We can use a non-trivial line bundle with Hermitian form rather than the trivial one, and we will do this in the next section.
\\\\From $\Omega$ we get the following sesquilinear form\footnote{We assume linearity in the second component.} on $L^2(M,d\mu):$ 
\begin{equation}\label{ses}
    \langle \Psi_1| \Psi_2\rangle:=\int_{M\times M}\overline{\Psi_1(y)}\Psi_2(x)\,\Omega(x,y)\,d\mu(x)\,d\mu(y)\;,
\end{equation}
which we assume is bounded (this is automatic if $\mu(M)<\infty$). The Riesz representation theorem gives an equivalence between bounded sesquilinear forms and bounded linear operators on Hilbert spaces, and the corresponding operator is given by
\begin{equation}
\hat{\Omega}\Psi(x)=\int_M \Psi(y)\Omega(y,x)\,d\mu(y)\;.
\end{equation}
That is, given a sesquilinear form $A(\cdot,\cdot)$ on a Hilbert space $(\mathcal{H},\langle\cdot,\cdot\rangle)\,,$ there is a unique linear operator $Q$ such that
\begin{equation}
    A(\Psi_1,\Psi_2)=\langle \Psi_1, Q\Psi_2\rangle\;.
\end{equation}
The conditions assumed on $\Omega$ imply that $\hat{\Omega}$ is an orthogonal projection, ie. $\hat{\Omega}^2=\hat{\Omega}$ and $\hat{\Omega}^*=\hat{\Omega}\,,$ from which it follows that 
\begin{equation}
    \langle \Psi_1| \Psi_2\rangle\ge 0\;.
\end{equation}
\begin{definition}
The space of physical states (or quantum states), denoted $\mathcal{H}_{\textup{phys}}\,,$ is the $\lambda=1$ eigenspace of $\hat{\Omega}.$
\end{definition}
On $\mathcal{H}_{\textup{phys}}\,,$ \cref{ses} agrees with the $L^2(M,d\mu)$-inner product. This eigenspace is analogous to a space of K\"{a}hler-polarized sections of a prequantum line bundle, and we can identify it with $L^2(M,d\mu)/\sim\,,$ where $\Psi_1\sim\Psi_2$ if $\hat{\Omega}(\Psi_1)=\hat{\Omega}(\Psi_2)\,.$ 
\\\\Such $\Omega$ are sometimes called reproducing kernels and have been widely studied, eg. \cite{sawin}, \cite{ma}. $\mathcal{H}_{\textup{phys}}$ has the special property that pointwise evaluation is a bounded linear functional, since by the definition of $\hat{\Omega},$ for $\Psi\in \mathcal{H}_{\textup{phys}}$
\begin{equation}
    |\Psi(x)|\le \|\Omega\|\|\Psi\|_{L^2}
\end{equation}
for all $x\in M.$ Therefore, by the Riesz representation theorem, $\mathcal{H}_{\textup{phys}}$ has ``delta functions" for each point $x\in M,$ denoted $|x\rangle.$ That is, if $\Psi\in\mathcal{H}_{\textup{phys}}\,,$ then there are genuine states $|x\rangle \in\mathcal{H}_{\textup{phys}}$ such that
\begin{equation}
    \langle x|\Psi\rangle=\Psi(x)\;.
\end{equation}
\begin{definition}
The states $\{|x\rangle\}_{x\in M}$ are the coherent states
\end{definition}
Coherent states form an overcomplete basis of $\mathcal{H}_{\textup{phys}}\,,$ in the sense that they are linearly dependent.  
\\\\It follows from the definition of $|x\rangle\,,$ \cref{ses} and the fact that for each $y\in M\,,$ $\Omega(x,y)\in \mathcal{H}_{\textup{phys}}\,,$ that
\begin{equation}
    \langle y| x\rangle=\Omega(x,y)\;,
\end{equation}
where $\langle y|$ is the dual linear functional of $|y\rangle\,.$ 
\begin{definition}
 There is a quantization map 
 \begin{equation}
  L^{\infty}(M,d\mu)\to \mathcal{B}(L^2(M,d\mu))\;,\;\;f\mapsto\hat{f}   
 \end{equation} 
 given by associating to each $f$ the following sesquilinear form:
\begin{equation}
    (\Psi_1,\Psi_2)\mapsto \int_{M^3}\overline{\Psi_1(z)}\Psi_2(x)\,f(y)\,\Omega(x,y)\Omega(y,z)\,d\mu(x)\,d\mu(y)\,d\mu(z)\;.
\end{equation}   
\end{definition}
This is the ``matrix element" of the observable $\gamma\mapsto f(\gamma(t))$ determined by the path integral (for any $t\in (0,1)$) and is given by
\begin{equation}
\int_{M\times M}\Bigg[\int_{\gamma(0)=x}^{\gamma(1)=y}f(\gamma(t))\,e^{i\int_0^1\gamma^*\omega}\,\mathcal{D}\gamma\Bigg]\,\overline{\Psi_1(y)}\Psi_2(x)\,d\mu(y)\,d\mu(x)\;.
\end{equation}
Explicitly,
\begin{equation}
    \hat{f}\Psi(x)=\int_{M^2}\Psi(z)\,f(y)\,\Omega(z,y)\Omega(y,x)\,d\mu(z)\,d\mu(y)\;.
\end{equation}
To be clear, the correspondence between the path integral and the operator formulations is an instance of the correspondence  determined by the Riesz representation theorem between sesquilinear forms and operators.\footnote{However, there is a sense in which the path integral formulation logically precedes the operator formalism. This is also one reason why the coherent space path integral (or phase space path integral, \cite{klein}) is nice, as this correspondence is obfuscated when a polarization is used.} 
\\\\These quantum operators have the property that $\hat{f}\hat{\Omega}=\hat{\Omega}\hat{f}=\hat{f}\,,$ hence
\begin{equation}
    \langle\Psi_1,\hat{f}\Psi_2\rangle=\langle \Psi_1|\hat{f}\Psi_2\rangle\;,
\end{equation}
where on the left is the inner product of $L^2(M,d\mu)\,.$ On $\mathcal{H}_{\textup{phys}}\,,$ we can write $\hat{f}$ diagonally in the coherent state basis, ie.
\begin{lemma}
\begin{equation}
    \hat{f}\vert_{\mathcal{H}_{\textup{phys}}}=\int_M f(x)|x\rangle\langle x|\,d\mu(x)\;.
\end{equation}
\end{lemma}
From this, it follows that $\mathcal{H}_{\textup{phys}}$ forms an irreducible representation of the quantum operators, so indeed $\mathcal{H}_{\textup{phys}}$ is a good Hilbert space.
\subsection{Example: Conventional Quantum Mechanics}\label{exmp}
Consider the symplectic manifold $\textup{T}^*\mathbb{R}$ with coordinates $(p,q)$ and symplectic form $\frac{dp\wedge dq}{2\pi\hbar}\,.$
Let
\begin{equation}
    \Omega(p,q,p',q')=\exp{\Big[-\frac{(p'-p)^2+(q'-q)^2}{4\hbar}+i\frac{pq'-qp'}{2\hbar}}\Big]\;.
\end{equation}
$\Omega$ is a propagator integrating $i\omega=\frac{i}{\hbar}\frac{p\,dq-q\,dp}{2}\,,$ according to \cref{prop1}. This is equal to 
\begin{equation}
    \langle p',q'|p,q\rangle\;.
\end{equation}
The state $|p,q\rangle$ is identified with the corresponding quantum mechanical coherent state, ie. the eigenvector of the lowering operator\footnote{These are usually defined in the context of the harmonic oscillator.} for which the expectation values of the momentum and position operators are $p,q,$ respectively. The physical Hilbert space $\mathcal{H}_{\textup{phys}}$ is the Segal-Bargmann representation (sometimes called the holomorphic representation), \cite{klauder3}.
\\\\Many other examples come from K\"{a}hler manifolds, as we will discuss in \cref{surf}.
\subsection{Relation to Deformation Quantization and Special Observables}
\subsubsection{Special Observables}
There is a special class of functions for which
\begin{equation}\label{canon}
    f(x)=\langle x|\hat{f}|x\rangle\;,
\end{equation}
where the latter is equal to 
\begin{equation}
    \int_M f(y)|\Omega(x,y)|^2\,d\mu(y)\;.
\end{equation}
Such functions generalize the harmonic functions on phase space (in particular, the coordinate functions $p,q$). There is a sense in which these observables are the most canonically quantizable ones — their classical and quantum expectation values agree on all states, in the sense of $C^*$-algebras (\cite{williams}). That is, given a normalized radon measure $dm$ on $M$ (ie. a mixed classical state), we get a mixed quantum state given by
\begin{equation}
T\in \mathcal{B}(\mathcal{H}_{\textup{phys}})\mapsto \int_M \langle x|T|x\rangle\,dm(x)\;,
\end{equation}
and if $f$ satisfies \cref{canon} then
\begin{equation}
\int_M \langle x|\hat{f}|x\rangle\,dm(x)=\int_M f(x)\,dm(x)\;.
\end{equation}
On the left is a quantum expectation value and on the right is a classical expectation value.\footnote{These describe states with classical statistical uncertainty.} In \cite{klauder}, the dual quantization scheme is used, ie. a function $f$ quantizes to an operator $\hat{f}$ such that $f(x)=\langle x|\hat{f}|x\rangle\,.$ However, such a quantization map doesn't need to exist. 
\\\\In special cases our quantization map
\begin{equation}
    L^{\infty}(M,d\mu)\to\mathcal{B}(\mathcal{H}_{\textup{phys}})
\end{equation}
is an isomorphism of vector spaces, and in this case we get a noncommutative product on $L^{\infty}(M,d\mu)\,.$ In some other case,  there may be a special subspace of functions for which the quantization map is injective and algebraically closed, in which case one gets a noncommutative product on this special subspace.
\subsubsection{Deformation Quantization}\label{dq}
We can get a quantization of a symplectic manifold $(M,\omega)$ by choosing a propagator integrating the connection on the prequantum line bundle, as discussed in \cref{mg}. One can often show that the quantization map is \textit{perturbatively} an isomorphism. Here, we should consider propagators integrating 
\begin{equation}
    \frac{i}{\hbar}\nabla\;,\;\,\frac{1}{\hbar}\in\mathbb{N}\;.
\end{equation}
In this case we obtain a formal deformation quantization, see equation (2.10) of \cite{brane}. This is what occurs in Toeplitz quantization \cite{bord}.  See \cref{kont}, \cite{sergio} for more on this relation. To make sense of what it measn to be an isomorphism perturbatively, we give the following definition:
\begin{definition}\label{nd}(see \cite{eli2})
Let $(M,\Pi)$ be a Poisson manifold and let $A\subset[0,1]$ be a set containing $0$ as an accumulation point. For each $\hbar\in A,$ let $M_{\hbar}$ be a unital $C^*$-algebra such that $M_{0}=L^{\infty}(M)$ and let 
\begin{equation}
    Q_{\hbar}:C_c^{\infty}(M)\to M_{\hbar}
\end{equation}
be $^*$- linear\footnote{This compatibility with $^*$ can be relaxed.} such that its image generates $M_{\hbar}\,.$ Furthermore, assume that $Q_0$ is the inclusion map. We say that $Q_{\hbar}$ is a (non-perturbative) deformation quantization of $(M,\Pi)$ if there is a star product $\star_{\hbar}$ on $C^{\infty}(M)[[\hbar]]$ such that, for all $n\in\mathbb{N},$
\begin{equation}
    \frac{1}{\hbar^{n}}||Q_{\hbar}(f)Q_{\hbar}(g)-Q_{\hbar}(f\star_{\hbar}^n g)||_{\hbar}\xrightarrow[]{\hbar \to 0} 0\;,
\end{equation}
where $f\star_{\hbar}^n g$ is the component of the star product up to order $n.$
\end{definition}
In this definition, we don't have to assume that $f\star_{\hbar}g$ is a star product a priori, just that it is an asymptotic expansion in $\hbar.$
\subsection{Non-Uniqueness of Primitives}\label{non}
It is important to note that, while primitives in the sense of \cref{cond} are unique (this uses connectivity of $\mathbb{R}$), primitives in the sense of \cref{condition1} are not. Therefore, we are abusing notation by writing \cref{ome} — we should really consider the right side to be the set of all functions satisfying the four aforementioned conditions and write
\begin{equation}
\Omega(x,y)\in\int_{\gamma(0)=x}^{\gamma(1)=y}e^{i\int_0^1\gamma^*\omega}\,\mathcal{D}\gamma\;.
\end{equation}
This is a similar abuse of notation one commits when stating that $f$ is \textit{the} antiderivative of $df\,.$ 
\\\\We address the non-uniqueness more in the example of \cref{surf}, where we show that the different propagators are distinguished by their second order Taylor expansions, rather than the expected first order.
\subsection{The Derivative of the Measure on Paths}\label{derp}
Consider a complex measure $d\gamma$ on the space of (discontinuous) paths $[0,1]\to (M,d\mu)\,.$ The condition corresponding to \cref{condition2} says that the derivative of the Radon-Nikodym derivative of the 2-dimensional distributions of $d\gamma$ equals $i\omega$ — its 2-dimensional distribution being the complex measure on $M\times M$ given by 
\begin{equation}
    \gamma_{2}(A):= \int_{(\gamma(0), \gamma(1))\in A}d\gamma\;.
\end{equation}
That is, the integral above is over all paths such that $(\gamma(0),\gamma(1))\in A\,.$ 
\\\\The Radon-Nikyodym derivative is taken with respect to the product measure $d\mu\times d\mu\,.$ This defines a function 
\begin{equation}
    \Omega:M\times M\to\mathbb{C}\;,
    \end{equation}
and we require that
\begin{equation}
    d_y\log{\Omega}(x,y)\vert_{y=x}=i\omega\vert_x\;.
\end{equation}
Of course, we are making some differentiability assumptions on the complex measure to do this. We generalize this to path integrals from higher dimensional simplices in \cref{pathm}.
\subsection{Extending the Fundamental Theorem of Calculus to Path Measures}
The ultimate goal is to extend the fundamental theorem of calculus in this setting. The main issues to understand are:
\begin{enumerate}
    \item non-uniqueness of the propagators differentiating to $i\omega$ (see \cref{surf} for a discussion),
    \item the problem of determining convergence on the right side of \cref{approx} for $\Omega$ which don't satisfy \cref{condition1}.
\end{enumerate}
The non-uniqueness is due to the fact that, under path integrals, we can't really integrate differential forms. This is because to uniquely specify Riemann-type integrals under path integrals we need to specify higher order information, eg. classically we integrate objects such as $f\,dx\,,$ but under path integrals we integrate objects such as 
\begin{equation}\label{high}
    f\,dx+g\,dx^2\;.
\end{equation}
Here, $dx^2=(dx)^2$ can be considered as a symmetric tensor, or as a homogenous polynomial of degree 2. We could integrate such objects classically as well, but we don't bother to make it explicit because the higher order terms don't contribute to the integral. However, in path integrals they do contribute. Picking good higher order terms often seems necessary to get convergence. To be a bit more precise, for a partition $x_0,\ldots,x_n$ of the interval $[0,1],$ classically
\begin{equation}\label{zero}
    \sum_{i=0}^{n-1}\Delta x_i^2\xrightarrow[]{\Delta x_i\to 0}0\;.
\end{equation}
This is why the usual definition of the Riemann integral doesn't depend on the choice of Riemann sum, ie. why the choice of $g$ in \cref{high} is never specified. However, \cref{zero} isn't quite true under path integrals and therefore we need to be more precise than just specifying a 1-form. We discuss this in detail in \cref{new}.
\section{The General Case: The Propagator of Symplectic Manifolds}\label{mg}
More generally, we can replace $L^2(M,d\mu)$ with square-integrable sections of a complex line bundle with Hermitian connection, where $i\omega$ of the previous section corresponds to the connection, and everything essentially works the same. Of course, the most natural setting for this is a prequantizable symplectic manifold $(M,\omega),$ which comes with a natural measure. Even better are K\"{a}hler manifolds, since there are canonically defined propagators, as discussed in \cref{examples}.
\\\\In this setting, consider a line bundle with Hermitian connection
\begin{equation}
    (\mathcal{L},\nabla,\langle \cdot,\cdot\rangle)\to (M,d\mu)\;,\footnote{Here, $d\mu$ is a measure on $M,$ which we really only care about up to multiplication by a positive, continuous function.}
\end{equation}
In this case, the path integral we want to compute is 
\begin{equation}\label{kernel}
   \Omega(x,y)=\int_{\gamma(0)=x}^{\gamma(1)=y}\mathcal{D}\gamma\,e^{\int_{0}^{1}\gamma^*\nabla}\;,
\end{equation}
where $e^{\int_{0}^{1}\gamma^*\nabla}$ denotes parallel transport between the vector spaces over $x$ and $y,$ over the curve $\gamma.$ That is, \cref{kernel} must define a linear map between any two vector spaces of the line bundle. The kernel $\Omega$ we need to consider then is a section of
\begin{equation}\label{sym}
    \pi_1^*\mathcal{L}^*\otimes\pi_2^*\mathcal{L}\to M\times M\;,\
\end{equation}
where $\pi_1,\pi_2:M\times M\to M$ are the projections onto the first and second factor, respectively. \Cref{sym} defines a line bundle over $M\times M\,,$ and points in this line bundle are naturally identified with linear maps between the corresponding vector spaces of $\mathcal{L}\,.$
\\\\The conditions of \cref{prop1} are replaced by: for all $(x,y)\in M\times M\,,$
\begin{align}
   \label{condition11} &\int_M \Omega(x,z)\circ \Omega(z,y)\,d\mu_z=\Omega(x,y)\,,
    \\\label{condition22}& \nabla_y{\Omega}(x,y)\vert_{y=x}=0\,,
    \\\label{last}&\Omega(x,x)=1\,,
    \\&\Omega(x,y)^*=\Omega(y,x)\,.
\end{align}
The subscript $y$ in \cref{condition22} means that the covariant derivative is being taken in the second component; each $l_x\in\mathcal{L}$ defines a section of $\mathcal{L}\,,$ given by
\begin{equation}
    y\mapsto l_x\Omega(x,y)\in\mathcal{L}_y\;,
\end{equation}
where the subscript is used to denote the corresponding vector space. \Cref{condition22} says that the covariant derivative of this section at $y=x$ is zero.
\begin{definition}\label{prop}
If $\Omega$ satisfies the previous conditions and the corresponding sesquilinear form \cref{ses2} is bounded, we say that it integrates $\nabla$ (as a propagator).\footnote{We may relax \cref{condition11}, see \cref{geo}.}
\end{definition}
Any section $\Omega$ of \cref{sym} which equals $1$ on the diagonal determines a connection on $\mathcal{L}$ (by differentiating in the second component at the diagonal),\footnote{We are using the identification of splittings $T\mathcal{L}\to\pi^*TM$ with connections.} and \cref{condition22} says that this connection is equal to $\nabla.$ Furthermore, for any smooth path $\gamma:[0,1]\to M$ and triangulation $0=t_0<\cdots<t_n=1$ of $[0,1],$ it follows that
\begin{equation}
    \prod_{i=0}^{n-1}\Omega(\gamma(t_i),\gamma(t_{i+1}))\;\;\footnote{This is essentially a generalized Riemann ``product" rather than a generalized Riemann sum, since a multiplicative structure is being used. See \cref{new}.}
\end{equation}
converges to parallel transport over $\gamma$ as $\Delta t_i\to 0.$ We discuss uniqueness of such a propagator in \cref{surf}.
\begin{remark}\label{symplectic}
In the symplectic case \cref{pcon}, we can interpret the exponent as specifying the connection only to leading order in $\hbar,$ in which case it is natural to relax \cref{condition22} to 
\begin{equation}
    \nabla_y\Omega(x,y)\vert_{y=x}=\mathcal{O}(\hbar)\;.
    \end{equation}
See \cref{pertu}.
\end{remark}
The next lemma tells us that that if we can find such an $\Omega$ without assuming \cref{last}, then we can always normalize so that all conditions are satisfied:
\begin{lemma}Suppose $(\nabla',\langle\cdot,\cdot\rangle',\Omega',d\mu')$ satisfies the following conditions: for all $(x,y)\in M\times M\,,$
\begin{align}
   \label{condition111} &\int_M \Omega'(x,z)\circ \Omega'(z,y)\,d\mu'_z=\Omega'(x,y)\,,
    \\\label{condition221}& \nabla'_y{\Omega'}(x,y)\vert_{y=x}=0\,,
    \\\label{last1}&\Omega'(x,x)\ne 0\,,
    \\&\Omega'(x,y)^*=\Omega'(y,x)\,.
\end{align}
Then $\Omega$ is a propagator integrating $\nabla,$ where $(\nabla, \langle\cdot,\cdot\rangle, \Omega,d\mu)$ is given by the following: let $f(x)=\Omega(x,x),$ then
\begin{align}
&d\mu=fd\mu'
\\&\Omega(x,y)=\frac{\Omega'(x,y)}{\sqrt{f(x)f(y)}}\,,
    \\& \nabla=\nabla'+\textup{dlog}\sqrt{f}\,,
    \\& \langle\cdot,\cdot\rangle=f\langle\cdot,\cdot\rangle'\;.
  \end{align}
\end{lemma}
The rest goes as in \cref{simple}, with the sesquilinear form \cref{ses} replaced by 
\begin{equation}\label{ses2}
\langle\Psi_1|\Psi_2\rangle=\int_{M\times M}\langle \Psi_1(y),\Psi_2(x)\Omega(x,y)\rangle\,d\mu_x\,d\mu_y\;.
\end{equation}
The only difference is that, due to nontriviality of the line bundle, the states $|x\rangle$ are only defined up to a scalar. To define a coherent state, choose a normalized vector $l_x$ over $x\,.$ We get a bounded linear functional on the physical Hilbert space $\mathcal{H}_{\textup{phys}}\,,$ which is given by square-integrable sections of the line bundle satisfying 
\begin{equation}\label{ac}
    \hat{\Omega}\Psi=\Psi\,,
    \end{equation}
where $\hat{\Omega}$ is the operator corresponding to the sesquilinear form \cref{ses2}. Given such a section $\Psi\,,$ 
\begin{equation}
    \Psi(x)=\lambda_\Psi\, l_x
    \end{equation}
for some $\lambda_{\Psi}=\mathbb{C}\,.$ The map
\begin{equation}\label{bd}
  \Psi\mapsto \lambda_{\Psi}  
\end{equation}
defines a bounded linear functional, by \cref{ac}. The Riesz representation theorem then defines the coherent states.
\begin{definition}
The coherent state $|l_x\rangle$ is defined to be the vector in $\mathcal{H}_{\textup{phys}}$ corresponding to the linear functional \cref{bd}.
\end{definition}
If $l_x=e^{i\theta}l'_x\,,$ then 
\begin{equation}
    |l_x\rangle\langle l_x|=|l'_x\rangle\langle l'_x|\;.
\end{equation}
Therefore, associated to $x$ there is a canonical projection operator, which we denote by $|x\rangle\langle x|\,.$\footnote{As pointed out to me by Francis Bischoff, the map $y\mapsto |y\rangle\langle y|$ is a version of the Kodaira embedding.} In addition, it follows that 
\begin{equation}\label{proa}
    l_x\Omega(x,y)= l_y\langle l_y|l_x\rangle\;,
\end{equation}
where $l_y$ is any vector over $y$ with unit norm.
\\\\Note that, it is only as a vector state that the coherent state  $|x\rangle$ isn't well-defined, as the linear functional on operators (ie. a $C^*$-algebra state)
\begin{equation}
    T\mapsto \langle x|T|x\rangle
\end{equation}
it is well-defined. Of course, the associated projection operator is also well-defined as discussed above.
\\\\Interestingly, if we can find a subspace of sections of the prequantum line bundle for which pointwise evaluation is continuous, then we can try to flip this construction to determine a propagator via \cref{proa}. We will do this in \cref{toe}.
\begin{remark}
Usually, a quantization is said to be something like an operator assignment $f\mapsto \hat{f}$ on a Hilbert space $\mathcal{H}_{\hbar}$, such that 
\begin{equation}
[\hat{f},\hat{g}]=i\hbar\{f,g\}+\mathcal{O}(\hbar^2)\;
\end{equation}
(eg. \cite{rieffel}). However, it's not clear that any such a structure should count as a quantization, since it appears to completely ignore any classical–quantum state correspondence. It is important that for any classical state $\rho$\footnote{In the sense of $C^*$-algebras, ie. $\rho$ is a continuous, positive linear functional on the $C^*$-algebra.} one can choose quantum states $\rho_{\hbar}$ such that $\rho_{\hbar}\to\rho$ as $\hbar\to 0.$ More precisely, we should have that 
\begin{equation}
    \rho_{\hbar}(\hat{f})\xrightarrow[]{\hbar\to 0}\rho(f)\;.
\end{equation}
Otherwise, one can't recover classical dynamics, as Heisenberg's equations only limit to Hamilton's equations on such sequences of states. The coherent state path integral approach comes with such states.
\end{remark}
\section{A New Look at the Riemann Integral}\label{new}
\subsection{The One-Dimensional Case}
To make sense of our definition of the path integral, we need to consider a more mathematically natural version of the Riemann integral. In particular, we need a version which behaves well with respect to pullbacks, since the path integral involves pulling back forms over paths. To motivate the definition, consider the following:
\begin{proposition}\label{int}(\cite{Lackman1}, \cite{Lackman2})
Let $f:[0,1]\to\mathbb{R}$ be a smooth function and let $F:[0,1]\times[0,1]\to\mathbb{R}$ be a smooth function which vanishes on the diagonal,\footnote{That is, $F(x,x)=0$ for all $x.$} such that
\begin{equation}\label{condi}
    \partial_y F(x,y)\vert_{y=x}=f(x)\;.
\end{equation}
Then 
\begin{equation}\label{rss}
    \sum_{i=0}^{n-1}F(x_i,x_{i+1})\xrightarrow[]{\Delta x_i\to 0}\int_0^1 f\,dx\;.
\end{equation}
\end{proposition}
Stated more geometrically, \cref{condi} says that the exterior derivative of $F$ in the second factor, evaluated at the diagonal, is equal to $f\,dx.$
\\\\The usual choices of $F$ are $F(x,y)=f(x)(y-x),\; f(y)(y-x)\,.$ These result in the left-point and right-point Riemmann sums, respectively. The left side of \cref{rss} can be considered to be a generalized Riemann sum. In particular, if $dF=f\,dx$ then 
\begin{equation}
    (x,y)\mapsto F(y)-F(x)
\end{equation}
also satisfies the conditions, so this definition makes the fundamental theorem of calculus almost tautological. This notion of integration easily to give a notion of Riemann sums of differential forms on manifolds. This gives a notion of integration on manifolds which mirrors the notion of integration on coordinate space, ie. it involves triangulations and not partitions of unity or coordinates.
\subsection{Integration Over Non-Differentiable Paths}\label{nond}
Importantly, consider a smooth map $\gamma:[0.1]\to[0,1]\,.$ This induces a map 
\begin{equation}
    \gamma\times\gamma:[0,1]\times[0,1]\to [0,1]\times[0,1]\;,
\end{equation}
and if $F$ satisfies the conditions of \cref{int} with respect to $f\,dx\,,$ then $(\gamma\times\gamma)^*F$ satisfies the conditions of \cref{int} with respect to $\gamma^*(f\,dx)\,.$ Therefore, if $F$ can be used to approximate 
\begin{equation}
    \int_0^1 f\,dx\;,
\end{equation}
then  $(\gamma\times\gamma)^*F$ can be used to approximate
\begin{equation}
    \int_0^1 \gamma^*(f\,dx)\;.
\end{equation}
As stated in \cref{non}, while the integral of $f\,dx$ makes sense classically, it isn't well-defined under path integrals because \cref{int} doesn't hold there. This is related to the non-differentiable nature of paths in the path integral. The problem is that, different functions $F$ chosen in \cref{condi} lead to different results under the path integral. Therefore, we need to specify more data. In the 1-dimensional Lagrangian formulation of the path integral (or better, its Euclidean version with respect to the Wiener measure), it is enough to specify 
\begin{equation}
  \partial^2_y F(x,y)\vert_{y=x}=g(x)\;. 
\end{equation}
That is, we need to specify the second derivative of $F$ as well. This means that 
\begin{equation}\label{well2}
    \int_0^1 \gamma^*(f\,dx+g\,dx^2)
\end{equation}
is well-defined in Wiener space $L^2(C[0,1],\mu_W)\,,$ if for a path $\gamma:[0,1]\to \mathbb{R}$ we take it to mean 
\begin{equation}\label{rs}
    \lim\limits_{\Delta t_i\to 0}\sum_{i=0}^{n-1}F(\gamma(t_i),\gamma(t_{i+1}))\;,
\end{equation}
where $F$ vanishes on the diagonal and 
\begin{equation}
\partial_y F(x,y)\vert_{y=x}=f(x)\,,\;\, \partial^2_y F(x,y)\vert_{y=x}=g(x)\;.
\end{equation}
To state a precise result, we have the following well-known result:
\begin{proposition}\label{wien}\footnote{The Riemann sum on the left defines the  Stratonovich integral and the Riemann sum on the right defines the It\^{o} integral, \cite{bernt}. The former satisfies the fundamental theorem of calculus.}
With respect to the Wiener measure on continuous paths $\gamma\in C([0,1],\mathbb{R})\,,$ and for a smooth function $f:\mathbb{R}\to\mathbb{R}\,,$ 
\begin{align}
    \nonumber &\sum_{i=0}^{n-1}f\bigg(\frac{\gamma(t_i)+\gamma(t_{i+1})}{2}\bigg)\,(\gamma(t_{i+1})-\gamma(t_i))\,-\,\sum_{i=0}^{n-1}f(\gamma(t_i))\,(\gamma(t_{i+1})-\gamma(t_i)) 
    \xrightarrow[]{\Delta t_i\to 0}\,\frac{1}{2}\int_0^1 \frac{df}{dx}(\gamma(t))\,dt
    \end{align}
in $L^2(C[0,1[).$
\end{proposition}
The left sum is the midpoint rule and the right sum is the left-point rule. The summands differ by a term of order $d\gamma^2\,.$ Classically, this term doesn't contribute to the integral because for smooth paths $d\gamma^2\sim dt^2\,.$ However, paths in Wiener space are generically H\"{o}lder continuous of exponent less than $1/2$, and $d\gamma^2\sim dt\,.$ Since the paths are generically H\"{o}lder continous with exponent greater than $1/3$ (\cite{rick}, \cite{peter}) it follows that all terms of order $d\gamma^3$ and higher are negligible. 
\\\\Note that, even on a smooth path these left and right sums are not Riemann sums in the traditional sense. In order to be a Riemann sum, the summand would need to be multiplied by $\Delta t_i,$ not $\Delta \gamma(t_i).$ However, in the more general sense of \cref{int}, they are Riemann sums. The same issue arises in Feynman's construction of the path integral. The sums involved in the approximations are not traditional Riemann sums, they are Riemann sums in our more general sense.
\\\\We will come back to this discussion in \cref{surf}.
\begin{remark}
It would be interesting to know if there are forms which are integrable in this sense but which aren't Lebesgue integrable. For example, consider a differentiable function $f:[0,1]\to\mathbb{R}$ such that $df$ isn't Lebesgue integrable. It is still the case that 
\begin{equation}
    (x,y)\mapsto f(y)-f(x)
\end{equation}
satisfies the conditions required to be used for a generalized Riemann sum of $df\,,$ and that 
\begin{equation}
    \sum_{i=0}^{n-1} f(x_{i+1})-f(x_i)
\end{equation}
converges in the limit (of course, to $f(1)-f(0)$). If we choose another function $F$ satisfying the conditions to be used for a Riemann sum of $df$ and such that 
\begin{equation}
    \sum_{i=0}^{n-1}F(x_i,x_{i+1}) 
\end{equation}
converges, must it converge to $f(1)-f(0)\,?$
\end{remark}
\subsection{Differentiating Path Measures to Differential Forms}\label{pathm}
Label the vertices of the (topological) standard $n$-simplex $|\Delta^n|$ by $0,\ldots,n.$ Suppose we have a complex measure $dX$ on the space of maps $X:|\Delta|^n\to M\,,$ with respect to some sigma algebra for which maps of the form
\begin{equation}
    X\mapsto f(X(i))
\end{equation}
are measurable, for $i=0,\ldots,n$ and for any $f:M\to\mathbb{C}\,.$
Then we can consider the finite dimensional distributions given by
\begin{equation}\label{meas}
A\subset {M^{n+1}}\mapsto\int_{(X(0),X(1),\ldots,X(n))\in A} dX\;.
\end{equation}
We can take the Radon-Nikodym derivative of these finite dimensional distributions with respect to $d\mu^{\times n+1},$ and this gives us a function 
\begin{equation}
    \Omega:M^{n+1}\to\mathbb{C}\;.
\end{equation}
\begin{definition}\label{condpath}
The derivative of $dX$ is the $n$-form on $M$ given by
\begin{equation}
    VE_0(\log{\Omega})\;.
\end{equation}
If we write $\Omega=\Omega(m_0,\ldots,m_n)\,,$ then $VE_0$ takes the exterior derivative of  $\log{\Omega}$ in each of the components $m_1,\ldots,m_n$ at $m_1=\cdots m_n=m_0.$ See \cref{des}.
\end{definition}
Of course, in the previous definition we are assuming that the finite dimensional distributions are absolutely continuous with respect to $d\mu^{\times n+1}$ and that it is $n$-times differentiable at the diagonal.
\\\\We suggest that any complex measure $dX$ on the space of maps $X:|\Delta^n|\to M$ which formally equals
\begin{equation}
    \mathcal{D}X\,e^{i\int_{|\Delta^n|}X^*\omega}
\end{equation}
should have a derivative which is equal to $i\omega.$ 
\section{Examples}\label{examples}
\subsection{Toeplitz Quantization and Bergman Kernels}\label{toe}
Let $(M,\omega,I)$ be a prequantizable K\"{a}hler manifold, with prequantum line bundle 
\begin{equation}
(\mathcal{L}_{},\nabla_{},\langle\cdot,\cdot\rangle_{})\to (M,i\omega/\hbar)\;.
    \end{equation}
Then pointwise evaluation of K\"{a}hler-polarized sections is a continuous map, and as previously discussed each unit-norm vector $l_x$ over $x$ determines a linear functional given by $\Psi\mapsto \lambda,$ where 
\begin{equation}
    \Psi(x)=\lambda l_x\;.
\end{equation}
We denote by $|l_x\rangle$ the associated section determined by the Riesz representation theorem. Note that, $\langle l_x|l_x\rangle$ is independent of $l_x,$ and we denote it by $\|x\|^2.$ 
\\\\We can define $\Omega_{}$ by
\begin{equation}
    l_x\Omega_{}(x,y)=l_y\frac{\langle l_y|l_x\rangle}{\|x\|\|y\|}\;,
\end{equation}
where $l_x, l_y$ have unit norm (this is independent of $l_y$).
\\\\The corresponding operator $\hat{\Omega}$ (\cref{corr}) is the orthogonal projection of $L^2(M,\mathcal{L})$ onto the subspace of holomorphic sections. The quantization map determined by the path integral is equal to the Toeplitz quantization map \cite{bord}. Its corresponding integral kernel is given by 
\begin{equation}
    \int_M f(z)\,\Omega(x,z)\circ\Omega(z,y)\,\|z\|^2\,\omega_z^n\;,
\end{equation}
and this is polarized with respect to the K\"{a}hler polarization on $M^{-}\times M$\footnote{To be more clear, it is polarized with respect to a possibly different connection — it is this integral kernel multiplied by $\|x\|\|y\|$ that is polarized with respect to $\nabla.$} (consistent with \cite{eli}).
\\\\$\Omega$ is closely related to the Bergman kernel $B$ (\cite{shiff}), which is important in the quantization of K\"{a}hler manifolds (\cite{sergio}), and $\|z\|^2\omega_z^n$ is called the Bergman measure (or coherent measure in \cite{sawin}). The Bergman kernel can be described as follows: let $\{\Psi_{\alpha}\}_{\alpha}$ be an orthonormal basis for the space of holomorphic sections. Then $B$ is a section of $\pi_1^*\mathcal{L}^*\otimes\pi_2^*\mathcal{L}\to M\times M\,,$ and is given by
\begin{equation}
    B(x,y)=\sum_{\alpha}\Psi_{\alpha}(x)^*\otimes \Psi_{\alpha}(y)\;.
\end{equation}
We then define a propagator\footnote{Note that, $\pi_1^*\mathcal{L}^*\otimes\pi_2^*\mathcal{L}$ is canonically trivial over the diagonal, so we can identify $x\mapsto B(x,x)$ with a positive-valued function.}
\begin{equation}
    \Omega(x,y)=\frac{B(x,y)}{\sqrt{B(x,x)}\sqrt{B(y,y)}}\;.
\end{equation}
Note that, $B(x,x)=\|x\|^2.$
\begin{lemma}\label{berg}
$\Omega$ integrates $\nabla,$ with respect to the Bergman measure 
\begin{equation}
   \frac{B(x,x)}{\hbar^n}\omega^n_x\;, 
\end{equation}
if and only if $x\mapsto B(x,x)$ is constant.
\begin{proof}
We need to show that $\nabla_y\Omega(x,y)\vert_{y=x}=0$ if and only if $B$ is constant on the diagonal. For the only if direction, differentiating along vectors in the $I=-i$ eigenspace and using the product rule implies that $x\mapsto B(x,x)$ is a holomorphic function. However, it is also real-valued, so it must be constant.
\\\\The other direction follows from product rule: since $\Omega(x,x)=1$ and $B(x,x)$ is constant, for vectors $Z\in T_{\mathbb{C}}M$ we have that
\begin{equation}
0=\nabla_{x,Z} B(x,y)\vert_{x=y}+\nabla_{y,Z}B(x,y)\vert_{x=y}\;.
\end{equation}
Since $B$ is antiholomorphic in $x,$ it follows that $\nabla_{y,Z}B(x,y)\vert_{x=y}=0$ for $Z$ in the $I=i$ eigenspace (ie. of the form $X-iI(X))\,.$ Since it is also true for $Z$ in the $I=-i$ eigenspace (since $B(x,y)$ is holomorphic in $y$) the result follows.
\end{proof}
\end{lemma}
\begin{exmp}
 $x\mapsto B(x,x)$ is constant for $M=\mathbb{C}P^n\,,$ where $\omega$ is the Fubini-Study symplectic form (\cite{shiff}).
\end{exmp} 
We get a function $M\times M\times M\to \mathbb{C}$ given by
\begin{equation}\label{comp}
    (x,y,z)\mapsto \Omega(x,y)\circ\Omega(y,z)\circ\Omega(z,x)\;.
\end{equation}
This defines a function on composable arrows of the pair groupoid, and if $B(x,x)$ is constant, then applying $VE_0\log$ to \cref{comp} gives $i\omega/\hbar.$ This is closely related to deformation quantization. \Cref{comp} should be interpreted as
\begin{equation}
    \int_{\gamma:S^1\to M}\mathcal{D}\gamma\,e^{\int_0^1\gamma^*\nabla}\;,
\end{equation}
where the domain of integration is all maps $\gamma:S^1\to M$ with three marked points on $S^1$ mapping to $x,y,z$ (compare with equations (2.9), (2.11) of \cite{brane}). If $M$ is simply connected then
\begin{equation}
    e^{\int_0^1\gamma^*\nabla}=e^{\frac{i}{\hbar}\int_D X^*\omega}\;,
\end{equation}
where $X$ is any map $X:D\to M$ which restricts to $\gamma$ on $\partial D.$ This corresponds to \cref{def} (after gauge reduction, see \cite{bon}).
\\\\For a relation of this quantization scheme to Berezin's star product and the Poisson sigma model on K\"{a}hler manifolds, see \cite{sergio}.
\\\\Note that, these examples somewhat contradicts the discussion in section 2.3 of \cite{brane}. This is due to the fact that the propagator isn't invariant under all symplectomorphisms, as is assumed.
\begin{remark}\label{pertu}
The map $x\mapsto B(x,x)$ is always asymptotically constant as $\hbar\to 0$ (\cite{sergio}), and $\Omega$ integrates a perturbed connection — the curvature of this perturbed connection is the Bergman K\"{a}hler form, which is an $\mathcal{O}(\hbar)$-perturbation of the original symplectic form, \cite{lu1}, \cite{ruan}. Therefore, it's still true that 
\begin{equation}
    \hbar\,VE_0\log{\big(\Omega(x,y)\circ\Omega(y,z)\circ\Omega(y,z)\big)}=i\omega+\mathcal{O}(\hbar)\;,
\end{equation}
which is what is needed for a quantization. See \cref{symplectic}. When $x\mapsto B(x,x)$ is constant we have an exact equality.\footnote{The property of the Bergman kernel being constant was used in the context of K\"{a}hler quantization  in \cite{cahen}. For it to be constant, it is enough that the manifold and line bundle are homogeneous.}
\\\\The function $x\mapsto B(x,x)$ is a well-studied map, and it being constant has other interesting consequences, eg. it implies stability in GIT (see page 2 of \cite{lu}). Metrics such that this map is constant are called balanced, \cite{donald}. For a generalization of the Bergman kernel to symplectic manifolds, see section 4 of \cite{ma} and \cite{sawin}.
\end{remark}
\subsection{Surfaces of Constant Curvature}\label{surf}
Here we will recall some explicit examples of the coherent state path integral \cref{pcon} computed by Daubechies, Klauder, on K\"{a}hler manifolds in \cite{klauder}, \cite{klauder2}. We will show that they satisfy \cref{prop}. To arrive at these expressions, the aforementioned authors (essentially) introduce Brownian motion into \cref{pcon} via a massive term $M,$ and send $M\to 0.$ \footnote{For more explicit examples of Bergman kernels, see \cite{sawin}.} In all of the following, the kernels are constant along the diagonal. In this paper, we have argued that this is true in general: whenever the Bergman kernel is constant along the diagonal, it computes the coherent state path integral.
\begin{exmp}
The following locally describe some 2-dimensional K\"{a}hler manifolds and their propagators \cref{pcon} in Darboux coordinates, for which the connection on the prequantum line bundle $(\mathcal{L},\nabla)$ is $pdq$ and $g$ is the Riemannian metric (which has constant curvature).
\begin{enumerate}
    \item Flat space: 
    \begin{equation}
         g=dp^2+dq^2\;,
        \end{equation}
        \begin{align}
         \Omega(p,q,p',q')=\exp{\Big[-\frac{(p'-p)^2+(q'-q)^2}{4\hbar}+i\frac{(p'+p)(q'-q)}{2\hbar}}\Big]\;.
    \end{align}
    
    \item Sphere: 
    \begin{equation}
        g=\frac{dp^2}{1-p^2/\hbar}+(1-p^2/\hbar)dq^2\;,
       \end{equation}
       \begin{align}
            \Omega(p,q,p',q')=\bigg[&\frac{1}{2}\Big(1+\frac{p'}{\sqrt{\hbar}}\Big)^{1/2}\Big(1+\frac{p}{\sqrt{\hbar}}\Big)^{1/2}e^{i\frac{q'-q}{2\sqrt{\hbar}}}
        +\\\nonumber\\ &\nonumber\frac{1}{2}\Big(1-\frac{p'}{\sqrt{\hbar}}\Big)^{1/2}\Big(1-\frac{p}{\sqrt{\hbar}}\Big)^{1/2}e^{i\frac{q'-q}{2\sqrt{\hbar}}}\bigg]^2\;.
    \end{align}
    \item Hyperbolic plane: \footnote{The examples computed by Klauder depend on parameters $s,\beta,$ which we've set equal to $1.$}
    \begin{equation}
         g=\hbar\frac{dp^2}{p^2}+\frac{p^2}{\hbar}\,dq^2\;,
        \end{equation}
        \begin{align}
        \Omega(p,q,p',q')=\Big(\frac{\hbar}{p'p}\Big)\bigg[\frac{\sqrt{\hbar}}{2p'}+\frac{\sqrt{\hbar}}{2p}-\frac{i(q'-q)}{2\sqrt{\hbar}}\bigg]^{-2}\;.
        \end{align}
\end{enumerate}
\end{exmp}
A simple computation shows that, indeed, 
\begin{equation}
    \textup{d}\textup{log}\,\Omega(p,q,p',q')\vert_{(p',q')=(p,q)}=\frac{i}{\hbar}p\,dq
    \end{equation}
(where the exterior derivative is taken in $(p',q')).$ Therefore, $\Omega$ integrates $\frac{i}{\hbar}\nabla$ in the sense of \cref{prop}. This illustrates the non-uniqueness of the value of \cref{pcon}, since the complement of an arc connecting the north and south pole of the sphere (which has measure zero) is symplectomorphic to $\mathbb{R}^2$ with its standard symplectic form.
\\\\Let us emphasize the following: there is a unique solution to \cref{cond}, \ref{cond2}, while the corresponding solution to \cref{condition1}-\ref{cond55} (or \cref{condition11}-\ref{last}) is not unique, and this is for the same reason that the Riemann integral doesn't have a unique value under path integrals, as discussed in \cref{nond} — in order to produce a unique value, one needs to specify a second order term, as in \cref{well2}. This leads us into the following:
\subsubsection{The Second Order Taylor Expansion of the Propagator (Bergman Kernel)}
Since the first order Taylor expansion doesn't uniquely determine the propagator, more interesting then is the second order Taylor expansion in the variables $(p',q')$ of $\log{\Omega}(p,q,p',q'),$ at $(p',q')=(p,q).$ This does distinguish between the different propagators. To make sense of this, we use the corresponding Riemannian metric to split the short exact sequence associated to the jet bundles
\begin{equation}\label{jet}
    0\to\text{Sym}^n T^*M\to J^n(M)\to J^{n-1}(M)\to 0\;.
    \end{equation}
We do this because, in order to define higher order Taylor expansions on manifolds, we need to split this sequence. Given a Riemannian metric, the result is the usual Taylor expansion formula, with derivatives replaced by symmetrized covariant derivatives, and terms like $(p'-p)^n$ replaced by $dp^n.$
\\\\For the hyperbolic plane, the non-zero Christoffel symbols are given by
\begin{equation}
    \Gamma^0_{00}=-\frac{1}{p}\,,\;\Gamma^0_{11}=-\frac{p^3}{\hbar^2}\,,\;\Gamma^{1}_{10}=\Gamma^{1}_{01}=\frac{1}{p}\;.
\end{equation}
\begin{proposition}
    In all three cases, the second order Taylor expansions are given by
\begin{equation}
   \frac{i}{\hbar}p\,dq-\frac{g}{4\hbar}\;. 
\end{equation}
\end{proposition}
Geometrically, this is equal to
\begin{equation}
    \frac{i}{\hbar}\nabla-\frac{g}{4\hbar}:\pi^*TM\to T\mathcal{L}\;,
\end{equation}
where we identify $g(X,X)\in\mathbb{C}$ with the corresponding vertical vector induced by the $\mathbb{C}^*$-action. 
\\\\Classically (ie. over smooth paths $\gamma$), the higher order term $-g/4\hbar$ doesn't contribute to the integral (see \cref{new}), ie.
\begin{equation}\label{dont}
    \int_0^1 \gamma^*\Big(\frac{i}{\hbar}p\,dq-\frac{g}{4\hbar}\Big)=\int_0^1\gamma^*\frac{i}{\hbar}p\,dq\;.
\end{equation}
Therefore, it is more precise to say that these propagators are equal to 
\begin{equation}\label{pre}
\int_{\gamma(0)=m}^{\gamma(1)=m'}\mathcal{D}\gamma\,e^{\frac{i}{\hbar}\int_0^1\gamma^*(\nabla+i\frac{g}{4})}\;,
\end{equation}
with the understanding that \cref{dont} doesn't hold on non-smooth paths. It would be interesting to know if the second order Taylor expansions (or 2-jets) uniquely determine the propagator — these 2-jets are given by a higher order van Est map (\cite{Lackman2}).
\subsubsection{The Lattice Approach to the Path Integral}\label{lat}
Similarly, it would be interesting to know which path integral approximations given by the lattice approach of \cite{Lackman2} converge to the propagator. This involves choosing a section $F$ of $\pi_1^*\mathcal{L}^*\otimes\pi_2^*\mathcal{L}$ over a neighborhood $U$ of the diagonal in $M\times M$ such that: $F$ restricted to the diagonal is equal to $1,$ and its derivative in the second component, evaluated on the diagonal, agrees with the connection — this means that the connection determined by $F$ agrees with $\nabla.$\footnote{We can always use the the Bergman kernel itself for $F$, in which case the following lattice construction clearly works.}
\\\\For example, we could cover $(M,\omega)$ with small enough open sets $\{O_{i}\}_i$ such that there is a unique geodesic connecting any two points in any of these open sets — this determines a neighborhood $U$ of the diagonal in $M\times M,$ with $(x,y)\in U$ if $x,y\in O_i$ for some $O_i.$ We get a section of $\pi_1^*\mathcal{L}^*\otimes\pi_2^*\mathcal{L}$ over $U,$ with the correct derivative, given by 
\begin{equation}
  F(m,m')= \exp{\bigg[-\frac{d(m,m')^2}{4\hbar}\bigg]}\exp{\bigg[\int_{m}^{m'}\gamma^*\nabla\bigg]}\;,\footnote{Let us emphasize that the term involving $d(m,m')$ is not an extra term in any way, it's just part of a term used in a (generalized) Riemann sum.}
\end{equation}
where $d$ is the distance between the corresponding points and $e^{\int_{m}^{m'}\gamma_{}^*\nabla}$ denotes parallel transport over the geodesic connecting $m$ to $m'.$ For a map $\gamma:[0,1]\to M$ and for a triangulation $0=t_0<\ldots<t_n=1$ of $[0,1],$ we get an approximation to the integrand of \cref{pre} given by
\begin{equation}
\prod_{i=0}^{n-1}F(\gamma(t_i),\gamma(t_{i+1}))\;.\footnote{The logarithm of this is a generalized Riemann sum of $\int_0^1\gamma^*\nabla$ over $[0,1].$}
\end{equation}
We can then compute the path integral \cref{pre} as in Feynman's approach, where the measure $\mathcal{D}\gamma$ is approximated using $\omega^n$ and the normalization is determined by the fact that the restriction of the path integral to the diagonal is equal to the $1.$ That is, we can try to compute
\begin{equation}
    \lim\limits_{n\to \infty}\frac{1}{C_n}\int_{U^{n}}(\pi_1^*\omega^n\times\pi_2^*\omega^n)\prod_{i=0}^{n}F(m_i,m_{i+1})\,\;,
\end{equation}
where $m_0=m,m_{n+1}=m',$ and $C_n$ is determined by the fact that for each $n$ the approximation is equal to $1$ if $m=m'.$
\\\\Note that indeed, the derivative (or pushforward) of $F$ in the second component, evaluated at the diagonal, agrees with $\nabla.$ Classicaly, this is enough information to get a well-defined integral (see \cref{new}) — the point is that, this is all of the information that is determined by the formal path integral. We can use the metric to split the short exact sequence associated to the jet bundle, \cref{jet}, and using this splitting $F$ also has the correct second derivative per \cref{pre}, ie. $-g/2\hbar.$
\section{Relation to Other Quantization Schemes}\label{geo}
\subsection{Relation to Kontsevich's Work, the Poisson Sigma Model and Weinstein's Program}\label{kont}
In the context of symplectic manifolds, the propagator $\Omega$ is a section of the prequantum line bundle over the symplectic groupoid, and it is an idempotent with respect to the twisted convolution algebra (\cite{eli}). Therefore, \cref{prop} naturally generalizes to Poisson manifolds. \Cref{condition22} implies that 
\begin{equation}
VE_0\big(\log{\big(\Omega(g_1)\circ\Omega(g_2)\circ\Omega(g_2^{-1}g_1^{-1})\big)}\big)=\frac{i}{\hbar}\Pi\;,    
\end{equation}
where $\Pi$ is the Poisson bivector and
\begin{equation}\label{comp}
    (g_1,g_2)\mapsto \Omega(g_1)\circ\Omega(g_2)\circ\Omega(g_2^{-1}g_1^{-1})
\end{equation}
is a function defined on composable arrows. To understand this function, consider three marked points $0,1,\infty$ on $S^1.$ \Cref{comp} should be interpreted as 
\begin{equation}
    \int_{g_1,g_2}\mathcal{D}\gamma \,e^{\frac{i}{\hbar}\int_D X^*\Pi}\;,
\end{equation}
where the integral is over morphisms $\gamma:TS^1\to T^*M$ whose restrictions to \begin{equation}
[\infty,0]\,,\, [0,1]\,,\, [1,\infty]\subset S^1
\end{equation}
are algebroid-homotopic to  $g_1,\, g_2,\, g_2^{-1}g_1^{-1},$ respectively, and where $X:TD\to T^*M$ is any morphism whose restriction to $\partial D$ is $\gamma.$\footnote{The symplectic groupoid being prequantizable implies that $\int_{S^2}X^*\Pi\in 2\pi\hbar\mathbb{Z}$ for morphisms  $X:TS^2\to T^*M,$ hence this is well-defined.}
Therefore, on a Poisson manifold this construction gives a $C^*$-algebra and quantization map, as in Weinstein's program \cite{weinstein}, and $\Omega$ can be used to compute 
\begin{equation}
    (f\star g)(m)=\int_{\gamma(\infty)=m}f(\gamma(1))\,g(\gamma(0))\,e^{\frac{i}{\hbar}\int_D X^*\Pi}\,\mathcal{D}\gamma\;.
    \end{equation}
Here, the integral is over all contractible morphisms $\gamma:TS^1\to T^*M,$
and $X:TD\to T^*M$ is any morphism whose restriction to $\partial D$ is $\gamma.$ This is formally equivalent to the Poisson sigma model approach to Kontsevich's star product \cite{kontsevich}, see \cite{bon}. In the symplectic case, it is the same as \cref{def} (after gauge reduction). See \cref{dq}.
\subsection{Relation to Kostant-Souriau's Geometric Quantization}
In \cref{exmp} we discussed the relation of these propagators to geometric quantization, using K\"{a}hler polarizations (also see \cref{kont}). However, one can ask about real polarizations.\footnote{Of course, there are examples where the propagator can be indirectly obtained using real polarizations, as discussed in \cref{exmp}. In this section we discuss viewing real-polarized sections as the physical states associated to a propagator.}  In \cite{Lackman5} (also see \cite{Lackman4}) we discussed a function 
\begin{equation}
    \Omega:\textup{T}^*\mathbb{R}\times \textup{T}^*\mathbb{R}\to\mathbb{C}
\end{equation}
satisfying all of the conditions to be a propagator integrating $i\frac{p\,dq-q\,dp}{2\hbar}\,,$ except that rather than satisfying $\Omega\ast\Omega=\Omega,$ it satisfies the relaxed condition
\begin{equation}
    \Omega\ast\Omega\ast\Omega=\Omega\;.
\end{equation}
Therefore, it's not a projection, but it still leads to good quantizations. In particular, the corresponding $\lambda=1$ eigenspace consists of all (complex)-linearly polarized sections (ie. for \textit{all} linear polarizations), and on this subspace $\Omega$ does give the correct inner product — the fact that all (complex) linearly-polarized sections are in this eigenspace leads to natural unitary equivalences between all such quantizations. Furthermore, it leads to a non-perturbative deformation quantization, which simultaneously acts on \textit{all} linearly-polarized sections.
\\\\It seems that we can view the failure of $\Omega\ast\Omega=\Omega$ as a failure of Fubini's theorem, since
\begin{equation}\label{propp}
    \Omega(p,q,p',q')=e^{i\frac{pq'-qp'}{2\hbar}}
\end{equation}
is not in $L^1(\textup{T}^*\mathbb{R})$ for fixed $(p,q).$ This $\Omega$ can still be considered to be a propagator, just not in the precise sense defined in this paper — though it seems this generalization is useful. It appears that any quantization scheme is (implicitly) 
attempting to assign a value to \cref{pcon}, though traditional geometric quantization forgets a lot of information about the quantization map.
\\\\\textbf{Question}: When does there exist an $\Omega$ integrating $i\nabla$ such that 
\begin{equation}
    \underbrace{\Omega\ast\Omega\ast\cdots\ast\Omega}_{n \text{ times}}=\Omega
\end{equation}
for some $n\ge 2\,?$
\begin{remark}
With respect to the standard K\"{a}hler metric, the Taylor expansion of the logarithm of \cref{propp} in the variables $(p',q'),$ at $(p',q')=(p,q),$ is exactly equal to 
\begin{equation}
    i\frac{p\,dq-q\,dp}{2\hbar}\;,
\end{equation}
ie. to all orders. In the language of \cite{Lackman2}, $VE_{\infty}(\log{\Omega})=p\,dq.$
\end{remark}
\begin{appendices}
\subsection{The van Est map and Riemann Sums on Manifolds}\label{Riemann}
The natural setting for the construction of the (generalized) Riemann sums on manifolds we give is the category of Lie groupoids, but it only involves the simplest Lie groupoid that exists, the pair groupoid. The pair groupoid $\textup{Pair}(M)\rightrightarrows M$ of a manifold $M$ is the groupoid with a unique arrow between any two points in $M.$ This is discussed more in \cite{Lackman1}, \cite{Lackman2}.
\begin{definition}\label{des}\footnote{This formulation is presented in \cite{Lackman2}, section 3. Also, \cite{Lackman1}. It is equivalent to, but a bit different from, the original formulation of Weinstein-Xu in \cite{weinstein1}. It was generalized  (to eg. $\mathbb{C}^*$-valued functions) in \cite{Lackman0}, which is relevant to path integrals.}
Let $M$ be a manifold. Let $\Omega:M^{n+1}\to\mathbb{C}$ be a smooth function which vanishes on the diagonal and which is invariant under even permutations. The van Est map (in degree 0) applied to $\Omega$ gives an $n$-form $VE_0(\Omega)$ on $M.$ Letting $X_1,\ldots X_n$ be vectors at a point $m\in M,$ it is defined by
\begin{equation}
    VE_0(\Omega)(X_1,\ldots,X_n)=n!\,X_n\cdots X_1\,\Omega(m,\cdot,\ldots,\cdot)\;,
\end{equation}
where $X_i$ differentiates $\Omega(m,\cdot,\ldots,\cdot)$ in the $ith$ component (where $m$ occupies the $0th$ component).
\end{definition}
\begin{exmp}
Let $M=\mathbb{R}^2$ and let $\Omega=\Omega(x_0,y_0,x_1,y_1,x_2,y_2)$ be a function on $\mathbb{R}^2\times\mathbb{R}^2\times\mathbb{R}^2$ which satisfies the conditions of the previous definition. Then
\begin{equation}
    VE_0(\Omega)(\partial_x\vert_{(x_0,y_0)},\partial_y\vert_{(x_0,y_0)})=2\,\partial_{y_2}\partial_{x_1}\Omega(x_0,y_0,\cdot,\cdot,\cdot,\cdot)\vert_{(x _1,y_1)=(x_2,y_2)=(x_0,y_0)}\;.
\end{equation}
\end{exmp}
States differently, $\Omega$ is a normalized $n$-cochain on $\textup{Pair}(M)$ which is invariant under even permutations.
\\\\Now, let $M$ be an oriented, $n$-dimensional manifold. Given a triangulation $\Delta_M$ of $M,$ up to even permutation there is a canonical ordering of the $(n+1)$ vertices of each $n$–dimensional simplex. Therefore, we can evaluate $\Omega:M^{n+1}\to\mathbb{C}$ on any $n$-dimensional face $\Delta\in\Delta_M$ by choosing such an ordering of the vertices and plugging them into $\Omega.$ This is well-defined because we are assuming that $\Omega$ is invariant under even permutations. We denote its value by $\Omega(\Delta).$
\begin{definition}(\cite{Lackman1}, \cite{Lackman2})
Let $M$ be an oriented $n$-dimensional manifold, let $\omega$ be an $n$-form on $M$ and let $VE_0(\Omega)=\omega.$ Then given a triangulation $\Delta_M$ of $M,$ the (generalized) Riemann sum of $\omega$ is defined to be
\begin{equation}
    \sum_{\Delta\in \Delta_M}\Omega(\Delta)\;,
\end{equation}
where the sum is over all $n$-dimensional simplices.
\end{definition}
\begin{theorem}\label{inte}(\cite{Lackman1}, \cite{Lackman2})
Suppose that $\Omega$ satisfies the conditions of \cref{des} and that $VE_0(\Omega)=\omega.$ Then
\begin{equation}
\sum_{\Delta \in\Delta_M }\Omega(\Delta)\xrightarrow[]{\Delta\to 0}\int_M\omega\;,
\end{equation}  
where the limit is taken over barycentric subdivisions of any triangulation $\Delta_M.$
\end{theorem}
The next proposition shows that this notion of Riemann sum is well–behaved with respect to pullbacks, which is important because pullbacks appear in path integrals:
\begin{proposition}
Let $\Delta_M$ be a triangulation of an $n$-dimensional manifold $M,$ let $f:M\to N$ be a smooth function and let $\omega$ be an $n$-form on $N$ with $VE_0(\Omega)=\omega.$ Then 
\begin{equation}
  \sum_{\Delta\in \Delta_M}f^*\Omega(\Delta)  
\end{equation}
is a (generalized) Riemann sum of $f^*\omega,$ where $f^*\Omega$ is the pullback via the map $M^{n+1}\to N^{n+1}$ induced by $f.$
\end{proposition}
\end{appendices}

 \end{document}